\pdfoutput=1
\RequirePackage{ifpdf}
\ifpdf 
\documentclass[pdftex]{sigma}
\else
\documentclass{sigma}
\fi
\numberwithin{equation}{section}

\newtheorem{Theorem}{Theorem}[section]
\newtheorem{Corollary}[Theorem]{Corollary}
\newtheorem{Lemma}[Theorem]{Lemma}
\newtheorem{Proposition}[Theorem]{Proposition}

\newcommand{\idem}{{\rm idem}}

\usepackage{bm}

\begin{document}

\allowdisplaybreaks

\newcommand{\arXivNumber}{1806.05375}

\renewcommand{\PaperNumber}{016}

\FirstPageHeading

\ShortArticleName{The $q$-Borel Sum of Divergent Basic Hypergeometric Series ${}_r\varphi_s(\bm{a};\bm{b};q,x)$}

\ArticleName{The $\boldsymbol{q}$-Borel Sum of Divergent Basic Hypergeometric\\ Series $\boldsymbol{{}_r\varphi_s(\bm{a};\bm{b};q,x)}$}

\Author{Shunya ADACHI}

\AuthorNameForHeading{S.~Adachi}

\Address{Graduate School of Education, Aichi University of Education, Kariya 448-8542, Japan}
\Email{\href{mailto:s217m064@auecc.aichi-edu.ac.jp}{s217m064@auecc.aichi-edu.ac.jp}, \href{mailto:s.adachi0324@icloud.com}{s.adachi0324@icloud.com}}

\ArticleDates{Received June 15, 2018, in final form February 24, 2019; Published online March 05, 2019}

\Abstract{We study the divergent basic hypergeometric series which is a $q$-analog of divergent hypergeometric series. This series formally satisfies the linear $q$-difference equation. In this paper, for that equation, we give an actual solution which admits basic hypergeometric series as a $q$-Gevrey asymptotic expansion. Such an actual solution is obtained by using $q$-Borel summability, which is a $q$-analog of Borel summability. Our result shows a $q$-analog of the Stokes phenomenon. Additionally, we show that letting $q\to1$ in our result gives the Borel sum of classical hypergeometric series.
The same problem was already considered by Dreyfus, but we note that our result is remarkably different from his one.}


\Classification{33D15; 39A13; 34M30}

\section{Introduction}
We set $q\in\mathbb{C}^*$ and $0<|q|<1$ throughout this paper. Basic hypergeometric series with the base $q$ is defined by
\begin{gather}\label{qHGS}
{}_r\varphi_s\left(\bm{a};\bm{b};q,x\right)=
{}_r\varphi_s\left(\begin{matrix}\bm{a}\\ \bm{b}\end{matrix};q,x\right)
:=\sum_{n\ge0} \frac{(a_1,a_2,\ldots,a_r;q)_n}{(b_1,\ldots,b_s;q)_n(q;q)_n}\big\{(-1)^nq^{\frac{n(n-1)}{2}}\big\}^{1+s-r}x^n,
\end{gather}
where $x\in\mathbb{C}$, $\bm{a}=(a_1,a_2,\ldots,a_r)\in\mathbb{C}^r$, $\bm{b}=(b_1,\ldots,b_s)\in\mathbb{C}^s$ and $b_1,\ldots,b_s\notin q^{-\mathbb{N}}$.

Basic hypergeometric series \eqref{qHGS} is a $q$-analog of the classical hypergeometric series
\begin{gather}\label{HGS}
{}_rF_s\left(\bm{\alpha};\bm{\beta};x\right)=
{}_rF_s\left(\begin{matrix}
\bm{\alpha}\\
\bm{\beta}
\end{matrix};x\right)
:=\sum_{n\ge0}\frac{(\alpha_1)_n(\alpha_2)_n\cdots(\alpha_r)_n}{(\beta_1)_n\cdots(\beta_s)_nn!}x^n,
\end{gather}
where $\bm{\alpha}=(\alpha_1,\alpha_2,\ldots,\alpha_r)\in\mathbb{C}^r$, $\bm{\beta}=(\beta_1,\ldots,\beta_s)\in\mathbb{C}^s$ and $\beta_1,\ldots,\beta_s\notin \mathbb{Z}_{\le0}$.

The radius of convergence of the series \eqref{qHGS} and \eqref{HGS} are both $\infty$, $1$ or $0$ according to whether $r<s+1$, $r=s+1$ or $r>s+1$. In this paper, we assume $r>s+1$ and $a_1a_2\cdots a_rb_1b_2\cdots b_s\neq0$. Therefore the series \eqref{qHGS} and \eqref{HGS} are divergent in this paper.

Basic hypergeometric series \eqref{qHGS} formally satisfies the following linear $q$-difference equation:
\begin{gather}\label{qHGEq}
\left[x(-\sigma_q)^{1+s-r}\prod_{j=1}^r(1-a_j\sigma_q)-(1-\sigma_q)\prod_{k=1}^s\left(1-\frac{b_k}{q}\sigma_q\right)\right]y(x)=0,
\end{gather}
where $\sigma_q$ is the $q$-shift operator defined by $\sigma_qy(x)=y(qx)$. As applying the local theory of linear $q$-difference equations (cf.\ Adams~\cite{Adams1929}), we see that the equation~\eqref{qHGEq} has a fundamental system of solutions around infinity\footnote{In his paper, Adams used the function $q^{-\frac{1}{2}(\log_q x)(\log_q x-1)}$ instead of $\theta_q(x)$. We remark that two functions $q^{-\frac{1}{2}(\log_q x)(\log_q x-1)}$ and $\theta_q(x)$ play the same role in constructing formal solutions since they satisfy the same $q$-difference equation
\begin{gather*}
\sigma_q^n y=x^{-n}q^{-\frac{n(n-1)}{2}}y.
\end{gather*}}:
\begin{gather}\label{qHGEq_Sol}
y_i^{(\infty)}(x)=\frac{\theta_q(-a_ix)}{\theta_q(-x)}{}_{r}\varphi_{r-1}
\left(\begin{matrix}
a_i,\frac{a_iq}{b_1},\frac{a_iq}{b_2},\ldots,\frac{a_iq}{b_s},0,\ldots,0\\
\frac{a_iq}{a_1},\ldots,\frac{a_iq}{a_{i-1}},\frac{a_iq}{a_{i+1}},\ldots,\frac{a_iq}{a_r}
\end{matrix};q,\frac{qb_1\cdots b_s}{a_1\cdots a_rx}\right),\qquad\!\! 1\le i\le r,\!\!\!
\end{gather}
which are convergent in $\vert x\vert>\vert qb_1\cdots b_s/a_1\cdots a_r\vert$.

In this paper, we construct an actual solution of \eqref{qHGEq} which admits \eqref{qHGS} as an asymptotic expansion and give its analytic continuation, by using $q$-Borel--Laplace transform (see Theorem~\ref{Theorem:Main1}). This result shows a $q$-analog of Stokes phenomenon. After that, we consider to take the limit $q\to1$ in Theorem~\ref{Theorem:Main1} (see Theorem~\ref{Theorem:Main2}).

The motivation of this study comes from as follows: Ichinobe \cite{Ichinobe2001} obtained the Borel sum of divergent classical hypergeometric series.
\begin{Theorem}[Ichinobe {\cite[Theorem 2.1]{Ichinobe2001}}]\label{Theorem:Ichinobe}
Let $\bm{\alpha}=(\alpha_1,\alpha_2,\ldots,\alpha_r)\in\mathbb{C}^r$ and $\bm{\beta}=(\beta_1,\beta_2,\ldots,$ $\beta_s)\in\mathbb{C}^s$. Assume $\alpha_i-\alpha_j\notin\mathbb{Z}$, $i\neq j$. Then ${}_rF_s(\bm{\alpha};\bm{\beta};x)$ is $1/(r-s-1)$-summable in any direction $d$ such that $d\neq0$ $({\rm mod}~2\pi)$ and its Borel sum $f(x)$ is given by
\begin{gather*}\label{Ichinobe_result1}
f(x)=C_{\bm{\alpha\beta}}\sum_{j=1}^rC_{\bm{\alpha\beta}}(j)(-x)^{-\alpha_j}{}_{s+1}F_{r-1}\left( \begin{matrix}
\alpha_j,1+\alpha_j-\bm{\beta}\\
1+\alpha_j-\widehat{\bm{\alpha}_j}
\end{matrix}
;\frac{(-1)^{1+s-r}}{x}\right),
\end{gather*}
where $x\in S(\pi,(r-s+1)\pi,\infty):=\{x\in\mathbb{C}^*;\left|\pi-\arg x\right|<(r-s+1)\pi/2)\}$ and $\widehat{\bm{\alpha}_j}\in\mathbb{C}^{r-1}$ is the vector which is obtained by omitting the $j$-th component from $\bm{\alpha}$ and
\begin{gather}\label{Ichinobe_result2}
C_{\bm{\alpha\beta}}=\frac{\Gamma(\bm{\beta})}{\Gamma(\bm{\alpha})},\qquad C_{\bm{\alpha\beta}}(j)=\frac{\Gamma(\alpha_j)\Gamma(\widehat{\bm{\alpha}_j}-\alpha_j)}{\Gamma(\bm{\beta}-\alpha_j)}.
\end{gather}
Here we use the following abbreviations
\begin{gather*}
\Gamma(\bm{\alpha})=\prod_{l=1}^r\Gamma(\alpha_l),\qquad \Gamma(\widehat{\bm{\alpha}_j}-\alpha_j)=\prod_{l=1,\,l\neq j}^r\Gamma(\alpha_l-\alpha_j).
\end{gather*}
\end{Theorem}
Since basic hypergeometric series \eqref{qHGS} is a $q$-analog of hypergeometric series \eqref{HGS}, whether there exists a $q$-analog of Theorem \ref{Theorem:Ichinobe} is a natural question. To answer this question, we use the theory of $q$-Borel summability. It is a $q$-analog of the theory of Borel summability and studied by many authors (cf.~Ramis \cite{Ramis1992}, Zhang \cite{Zhang2002,Zhang2005}, Di Vizio--Zhang \cite{Vizio-Zhang} and Dreyfus--Eloy~\cite{DreyfusEloy2016}). Dreyfus~\cite{Dreyfus2015_1} proved that for every formal power series solution of a linear $q$-difference equation with rational coefficients, we may construct a meromorphic solution of the same equation applying several $q$-Borel and $q$-Laplace transformations of appropriate orders and appropriate direction.

Our results are generalizations of Zhang \cite{Zhang2002} and Morita \cite{Morita2014}. Zhang studied the $q$-Borel summability of divergent basic hypergeometric series ${}_2\varphi_0$ which is the case of $r=2$ and $s=0$ in our result. Later, Morita obtained similar results with Zhang for ${}_3\varphi_1$, which is the case of $r=3$ and $s=1$ in our result. The common point of their assumptions is that $r$ and $s$ satisfy $r-s=2$. Our result gives the resummation of ${}_r\varphi_s$ in the case of $r-s\ge2$.

This paper is organized as follows. We fix our notions and review the theory of $q$-Borel summability in Section \ref{Section:Preliminary}. Main results of this paper are given in Section \ref{Section:Main Results}. In Sections~\ref{Section:Proof_Main1} and~\ref{Section:Proof_Lemma:Main1}, we give proofs of Theorem~\ref{Theorem:Main1} and Lemma~\ref{Lemma:Main1} respectively. A~proof of Theorem~\ref{Theorem:Main2} is given in Section~\ref{Section:Proof_Main2}.

\subsection*{On an article of Dreyfus}
We note that Dreyfus \cite{Dreyfus2015_2} treated the resummation of basic hypergeometric series \eqref{qHGS} via $q$-Borel summability. But his result is remarkably different from our result: In Dreyfus' paper, the order of $q$-Laplace transform and $q$-Borel transform are different. Therefore the function he obtained doesn't satisfy the same equation with \eqref{qHGS}. To explain this more detail, we use the following notations which are used in Dreyfus' paper \cite{Dreyfus2015_2} only in this subsection.
\begin{itemize}\itemsep=0pt
\item $q>1$ and $p=1/p$,
\item $\underline{q}:=q^{1/(r-s-1)}$ and $\underline{p}:=p^{1/(r-s-1)}$,
\item the definitions of $\hat{B}_q$ and $L_{q}^{[d]}$ are essentially same with this paper (see Section~\ref{Review}).
\end{itemize}
Under these notations, Dreyfus considered the resummation of ${}_r\varphi_s(\bm{a};\bm{b};\underline{p},z)$. He said that~\cite[pp.~478--479]{Dreyfus2015_2}
\begin{gather*}
\text{\dots, we can compute a solution of the same linear $\displaystyle\sigma_{\underline{q}}$-equation than ${}_r\varphi_s\!\left(
\begin{matrix}
a_1,\ldots,a_r\\
b_1,\ldots,b_s
\end{matrix};\underline{p},z
\right)$}\\
\text{$=\sum\limits_{n\in\mathbb{N}}\!\frac{(a_1,{\ldots},a_r;\underline{p})_n}{(b_1,{\ldots},b_s;\underline{p})_n(\underline{p};\underline{p})_n}\big\{(-1)^n\underline{p}^{\frac{n(n-1)}{2}}\big\}^{1+s-r}z^n$ applying successively to it $\hat{B}_{\underline{q}}$~and~$L_{\underline{q}}^{[d]}$.}
\end{gather*}
However, since
\begin{gather*}
\left(\underline{p}^{\frac{n(n-1)}{2}}\right)^{1+s-r}=p^{-\frac{n(n-1)}{2}}=q^{\frac{n(n-1)}{2}},
\end{gather*}
we have to apply $q$-Borel transform $\hat{B}_q$ to ${}_r\varphi_s(\bm{a};\bm{b};\underline{p},z)$, not $\underline{q}$-Borel transform $\hat{B}_{\underline{q}}$. In addition, the computation actually written in \cite{Dreyfus2015_2} seems to be wrong. He said that \cite[p.~479]{Dreyfus2015_2}
\begin{gather*}
\text{Applying $\hat{B}_{\underline{q}}$ to ${}_r\varphi_s\left(
\begin{matrix}
a_1,\ldots,a_r\\
b_1,\ldots,b_s
\end{matrix};\underline{p},z
\right)$, we obtain for all $d\not\equiv(r-s-1)\pi[2\pi]$}\\
h(\zeta):={}_r\varphi_{r-1}\left(\begin{matrix}
a_1,\ldots,a_r\\
b_1,\ldots,b_s,0,\ldots,0
\end{matrix};\underline{p},(-1)^{1+s-r}\zeta
\right)\in\mathbb{H}_{\underline{q},1}^d.
\end{gather*}
But actually the function $h(\zeta)$ is equal to $\hat{B}_q ({}_r\varphi_s)$, not $\hat{B}_{\underline{q}} ({}_r\varphi_s)$ (see Section \ref{Section:Proof_Main1}). Here we remark that the series $\hat{B}_{\underline{q}} ({}_r\varphi_s)$ must be still divergent because of $q>\underline{q}$. After that, the $\underline{q}$-Laplace transform of $h(\zeta)$ was calculated in Lemma~7.2. As a result, the resummation denoted by~$\mathbb{S}_{\underline{q}}^{[d]}({}_r\varphi_s)$ in Theorem~7.3 is actually equal to
\begin{gather}\label{Dreyfus_wrong}
\big(L_{\underline{q}}^{[d]}\circ\hat{B}_{q}\big)({}_r\varphi_s).
\end{gather}
Since the orders of $q$-Borel and $q$-Laplace transformations are different, the function~\eqref{Dreyfus_wrong} doesn't satisfy the same equation with ${}_r\varphi_s$.

However, we remark that the limit $q\to1$ of $\mathbb{S}_{\underline{q}}^{[d]}({}_r\varphi_s)$ in Theorem~7.3 is accidentally same with the Borel sum of (classical) hypergeometric series~\eqref{HGS}. In other words, the statement of Theorem~7.4 itself is correct, but the proof is incorrect.

In this paper we give the correct result about the resummation of basic hypergeometric series~\eqref{qHGS} in Theorem~\ref{Theorem:Main1}. In addition, we consider the relation between the resummation and local solution around infinity~\eqref{qHGEq_Sol}. We remark that this topic wasn't treated by Dreyfus. In Theorem~\ref{Theorem:Main2}, we give the limit $q\to1$ of the resummation of \eqref{qHGS} in Theorem~\ref{Theorem:Main1}, with a correct proof.

\section{Preliminary}\label{Section:Preliminary}
\subsection{Basic notions}
Let us fix our notations. We denote by $\bm{0}_n$ the $n$-vector $(0,0,\ldots,0)$.

For $n=0,1,2,\ldots,$ we set the Pochhammer symbol as
\begin{gather*}
(a)_n:=
\begin{cases}
1,& n=0,\\
a(a+1)\cdots(a+n-1), &n\geq1,
\end{cases}
\end{gather*}
and the $q$-shift factorial is defined by
\begin{gather*}
(a;q)_n:=\prod_{j=0}^{n-1}\big(1-aq^j\big),\qquad (a;q)_\infty:=\prod_{j=0}^\infty \big(1-aq^j\big).
\end{gather*}
For any $a\in\mathbb{C}$, the infinite product $(a;q)_\infty$ is convergent. We set
\begin{gather*} (a_1,a_2,\ldots,a_r;q)_{n}:=\prod\limits_{j=1}^r (a_j;q)_n \qquad \text{for $n=0,1,2,\ldots$ or $n=\infty$}.
\end{gather*}

The theta function with the base $q$ is defined by
\begin{gather*}
\theta_q(x):=\sum_{n\in\mathbb{Z}}q^{\frac{n(n-1)}{2}}x^n,
\end{gather*}
which is holomorphic on $\mathbb{C}^*=\mathbb{C}\setminus\{0\}$. The following properties hold:
\begin{gather}
\theta_q(x) =(q,-x,-q/x;q)_\infty,\label{theta_tpf}\\
\theta_q(q^nx) =x^{-n}q^{-\frac{n(n-1)}{2}}\theta_q(x),\qquad n\in\mathbb{Z},\label{theta_q-diff}\\
\theta_q(qx) =\frac{\theta_q(x)}{x}=\theta_q\left(\frac{1}{x}\right)\label{theta_inverse}.
\end{gather}
For $\lambda\in\mathbb{C}^*$, we set the $q$-spiral $[\lambda;q]:=\lambda q^\mathbb{Z}=\{\lambda q^n ;\,n\in\mathbb{Z}\}$. From the equality~\eqref{theta_tpf}, we see
\begin{gather*}
\theta_q\left(-\frac{\lambda}{x}\right)=0\quad \Leftrightarrow \quad x\in[\lambda;q].
\end{gather*}
We remark that the $q$-spiral $[\lambda;q]$ is identified with an element of $\mathbb{C}^*/q^\mathbb{Z}$.

The $q$-gamma function $\Gamma_q(x)$ is defined by
\begin{gather*}
\Gamma_q(x):=\frac{(q;q)_\infty}{(q^x;q)_\infty}(1-q)^{1-x}.
\end{gather*}
When $q$ is a real number and satisfies $0<q<1$, the limit of $\Gamma_q(x)$ as $q\to1$ gives the gamma function (cf.\ Gasper--Rahman {\cite[Section~1.10]{Gasper&Rahman2004}})
\begin{gather}\label{qgamma_reduction}
\lim_{q\to1}\Gamma_q(x)=\Gamma(x).
\end{gather}

\subsection[Review of $q$-Borel summability]{Review of $\boldsymbol{q}$-Borel summability}\label{Review}
We review the theory of $q$-Borel summability (cf.~Zhang \cite{Zhang2002, Zhang2005} and Dreyfus--Eloy~\cite{DreyfusEloy2016}).

Let $\mathbb{C}[[x]]$ be the ring of formal power series of $x$ and $\hat{f}(x)=\sum\limits_{n\ge0}a_nx^n\in\mathbb{C}[[x]]$. For $\hat{f}(x)$, we define its $q$-Borel transform $\hat{\mathcal{B}}_{q;1}\colon \mathbb{C}[[x]]\to\mathbb{C}[[\xi]]$ by
\begin{gather*}
(\hat{\mathcal{B}}_{q;1}\hat{f})(\xi):=\sum_{n\ge0}a_nq^{\frac{n(n-1)}{2}}\xi^n.
\end{gather*}

We denote $\mathcal{M}(\mathbb{C}^*,0)$ the field of functions that are meromorphic on some punctured neigh\-borhood of~0 in $\mathbb{C}^*$. More explicitly $f\in\mathcal{M}(\mathbb{C}^*,0)$ if and only if there exists $V$, an open neighborhood of 0, such that $f$ is analytic on $V\setminus\{0\}$.

Let $[\lambda;q]\in\mathbb{C}^*/q^\mathbb{Z}$ and $f\in\mathcal{M}(\mathbb{C}^*,0)$. It is said that $f$ belongs to $\mathbb{H}_{q;1}^{[\lambda;q]}$ if there exist a positive constant $\varepsilon$ and a domain $\Omega\subset\mathbb{C}$ such that:
\begin{itemize}\itemsep=0pt
\item $\displaystyle\bigcup_{m\in\mathbb{Z}}\left\{x\in\mathbb{C}^* ;\,\left|x-\lambda q^m\right|<\varepsilon\left|q^m\lambda\right|\right\}\subset\Omega$.
\item The function $f$ can be continued to an analytic function on $\Omega$ with $q$-exponential growth at infinity, which means that there exist positive constants $L$ and $M$ such that for any $x\in\Omega\setminus\{0\}$, the following holds:
\begin{gather*}
|f(x)|<L\theta_{|q|}(M|x|).
\end{gather*}
We note that this estimate does not depend on the choice of $\lambda$.
\end{itemize}

Let $[\lambda;q]\in\mathbb{C}^*/q^{\mathbb{Z}}$. For $f\in\mathbb{H}_{q;1}^{[\lambda;q]}$, we define its $q$-Laplace transform $\mathcal{L}_{q;1}^{[\lambda;q]}\colon \mathbb{H}_{q;1}^{[\lambda;q]}\to\mathcal{M}(\mathbb{C}^*,0)$ by
\begin{gather*}
\big(\mathcal{L}_{q;1}^{[\lambda;q]}f\big)(x):=\frac{1}{1-q}\int_{0}^{\lambda\infty}\frac{f(\xi)}{\theta_q\big(\frac{\xi}{x}\big)}\frac{d_q\xi}{\xi}=\sum_{n\in\mathbb{Z}}\frac{f(\lambda q^n)}{\theta_q\big(\frac{\lambda q^n}{x}\big)}.
\end{gather*}
Here, this transformation is given in Jackson integral (cf.\ Gaspar--Rahman \cite[Section~1.10]{Gasper&Rahman2004}). For sufficient small $|x|$, the function $\big(\mathcal{L}_{q;1}^{[\lambda;q]}f\big)(x)$ has poles of order at most one that are contained in the $q$-spiral $[-\lambda;q]$. We remark that the well-definedness of $q$-Laplace transform can be shown in the same way as Th\'eor\`eme~1.3.2 of Zhang~\cite{Zhang2002}. If the function $\hat{f}(x)\in\mathbb{C}[[x]]$ satisfies $\big(\hat{\mathcal{B}}_{q;1}\hat{f}\big)\in\mathbb{H}_{q;1}^{[\lambda;q]}$, we call that $\hat{f}(x)$ is $[\lambda;q]$-summable and $\big(\mathcal{L}_{q;1}^{[\lambda;q]}\circ\hat{\mathcal{B}}_{q;1}\hat{f}\big)(x)$ is the $[\lambda;q]$-sum of~$\hat{f}(x)$.

We give some fundamental properties of $q$-Borel--Laplace transform without proofs.
\begin{Proposition}[Dreyfus--Eloy {\cite[Section~1]{DreyfusEloy2016}}] Let $[\lambda;q]\in\mathbb{C}^*/q^\mathbb{Z}$ and $f$ be an analytic function. Then, the followings hold:
\begin{itemize}\itemsep=0pt
\item $\big(\hat{\mathcal{B}}_{q;1}f\big)\in\mathbb{H}_{q;1}^{[\lambda;q]}$,
\item $\big(\mathcal{L}_{q;1}^{[\lambda;q]}\circ\hat{\mathcal{B}}_{q;1}f\big)=f$.
\end{itemize}
\end{Proposition}
\begin{Proposition}\label{prop_borel2} Let $[\lambda;q]\in\mathbb{C}/q^\mathbb{Z}$. For $\hat{f}(x)\in\mathbb{C}[[x]]$ and $g\in\mathbb{H}_{q;1}^{[\lambda;q]}$, the followings hold.
\begin{itemize}\itemsep=0pt
\item $\big(\hat{\mathcal{B}}_{q;1}\big(x^m\sigma_q^n\hat{f}\big)\big)=q^{\frac{m(m-1)}{2}}\xi^m\sigma_q^{n+m}\big(\hat{\mathcal{B}}_{q;1}\hat{f}\big)$,
\item $\big(\mathcal{L}_{q;1}^{[\lambda;q]}(\xi^m\sigma_q^ng)\big)=q^{-\frac{m(m-1)}{2}}x^m\sigma_q^{n-m}\big(\mathcal{L}_{q;1}^{[\lambda;q]}g\big)$.
\end{itemize}
\end{Proposition}
The first equality is shown by direct calculation. The second equality can be proved by using the property (2.2) of theta function. From this proposition, we have the following important corollary.
\begin{Corollary}\label{cor_q-diff}
Let $\hat{f}(x)\in\mathbb{C}[[x]]$ be a formal solution of a~linear $q$-difference equation. If $\hat{f}(x)$ is $[\lambda;q]$-summable, its $[\lambda;q]$-sum satisfies the same $q$-difference equation with $\hat{f}(x)$.
\end{Corollary}

At the end of this section, we give asymptotic properties of the $q$-Borel--Laplace transform.

Let $[\lambda;q]\in\mathbb{C}^*/q^\mathbb{Z}$, $\hat{f}(x)=\sum\limits_{n\ge0}a_nx^n\in\mathbb{C}[[x]]$ and $f\in\mathcal{M}(\mathbb{C}^*,0)$. Then we define
\begin{gather*}
f\sim_1^{[\lambda;q]}\hat{f},
\end{gather*}
if for any positive numbers $\varepsilon$ and $R$, there exist positive constants $C$ and $K$ such that for any $N\ge1$ and
\begin{gather*}
x\in\left\{x\in\mathbb{C}^* ;\,|x|<R\right\}\setminus\bigcup_{m\in\mathbb{Z}}\left\{x\in\mathbb{C}^* ;\,|x+\lambda q^m|<\varepsilon |q^m\lambda |\right\},
\end{gather*}
we have
\begin{gather*}
\left|f(x)-\sum_{n=0}^{N-1}a_nx^n\right|\le LM^N|q|^{-\frac{N(N-1)}{2}}|x|^N.
\end{gather*}
\begin{Proposition}[Dreyfus--Eloy {\cite[Section~1]{DreyfusEloy2016}}]\label{qBorel_AE}
Let $[\lambda;q]\in\mathbb{C}^*/q^\mathbb{Z}$ and $\hat{f}\in\mathbb{C}[[x]]$ with $\big(\hat{\mathcal{B}}_{q;1}\hat{f}\big)\in\mathbb{H}_{q;1}^{[\lambda;q]}$. Then $\big(\mathcal{L}_{q;1}^{[\lambda;q]}\circ\hat{\mathcal{B}}_{q;1}\hat{f}\big)\sim_{1}^{[\lambda;q]}\hat{f}$.
\end{Proposition}

\section{Main results}\label{Section:Main Results}
In this section, we give our main results which are about the resummation of the basic hypergeometric series
\begin{gather}\label{qHGS2}
\hat{f}(x)={}_r\varphi_s(\bm{a};\bm{b};q,x)={}_r\varphi_s\left( \begin{matrix}\bm{a}\\ \bm{b} \end{matrix};q,x\right).
\end{gather}
We remark that $\hat{f}(x)$ is the formal power series solution of the linear $q$-difference equation \eqref{qHGEq}. Our results are stated as follows.
\begin{Theorem}\label{Theorem:Main1}
We put $k=r-s-1$ and $p=q^k$. Assuming $a_i/a_j\notin q^{\mathbb{Z}}$, $i\neq j$, basic hypergeometric series \eqref{qHGS2} is $[\lambda;p]$-summable for any $\lambda\in\mathbb{C}^*\setminus\big[(-1)^{k};q\big]$ and its $[\lambda;p]$-sum ${}_rf_s(\bm{a};\bm{b};\lambda;q,x)=\big(\mathcal{L}_{p;1}^{[\lambda:p]}\circ\mathcal{\hat{B}}_{p;1}\hat{f}\big)(x)$ is given by
\begin{gather}
{}_rf_s(\bm{a};\bm{b};\lambda;q,x) =\frac{(a_2,\ldots,a_r,b_1/a_1,\ldots,b_{s}/a_1;q)_\infty}{(b_1,\ldots,b_{s},a_2/a_1,\ldots,a_r/a_1;q)_\infty}\frac{\theta_{p}\big(pa_1^kx/\lambda\big)}{\theta_{p}(px/\lambda)}
\frac{\theta_q\big((-1)^{1-k}a_1\lambda\big)}{\theta_q\big((-1)^{1-k}\lambda\big)}\nonumber\\
\hphantom{{}_rf_s(\bm{a};\bm{b};\lambda;q,x) =}{}
\times{}_r\varphi_{r-1}\left(\begin{matrix}a_1,a_1q/b_1,\ldots,a_1q/b_s,\bm{0}_k\\a_1q/a_2,\ldots,a_1q/a_r
\end{matrix} ;q,\frac{qb_1\cdots b_s}{a_1\cdots a_rx}\right)\nonumber\\
\hphantom{{}_rf_s(\bm{a};\bm{b};\lambda;q,x) =}{}
+\idem(a_1;a_2,\ldots,a_r),\label{Main1}
\end{gather}
where $x\in\mathbb{C}^*\setminus[-\lambda;p]$ and $|qb_1b_2\cdots b_s/a_1\cdots a_rx|<1$. Here the symbol ``$\,\idem(a_1;a_2,\ldots,a_r)$'' means the sum of $r-1$ terms which are obtained by interchanging $a_1$ with each $a_2,\ldots,a_r$ in the preceding expression.
The function ${}_rf_s(\bm{a};\bm{b};\lambda;q,x)$ has simple poles on the $p$-spiral $[-\lambda;p]$.
\end{Theorem}
From Corollary~\ref{cor_q-diff} and the well-definedness of $q$-Laplace transform, we see that the $[\lambda;p]$-sum ${}_rf_s(\bm{a};\bm{b};\lambda;q,x)=\big(\mathcal{L}_{p;1}^{[\lambda:p]}\circ\mathcal{\hat{B}}_{p;1}\hat{f}\big)(x)$ itself becomes the meromorphic solution of~\eqref{qHGEq} on some punctured neighborhood of the origin in $\mathbb{C}^*$. The expression~\eqref{Main1} gives the analytic continuation of ${}_rf_s(\bm{a};\bm{b};\lambda;q,x)$ to the neighborhood of infinity.
 In addition, this shows a $q$-analog of the Stokes phenomenon. Indeed, an actual solution~\eqref{Main1} has~\eqref{qHGS2} as its asymptotic expansion (see Proposition~\ref{qBorel_AE}) and its value changes depending on the choice of $[\lambda;p]\in\mathbb{C}^*/p^{\mathbb{Z}}$.
Using the fundamental system of solutions around infinity \eqref{qHGEq_Sol}, the $[\lambda;p]$-sum \eqref{Main1} can be expressed as
\begin{gather*}
{}_rf_s(\bm{a};\bm{b};\lambda;q,x)=\sum_{j=1}^r M_j(x,\lambda)y_j^{(\infty)}(x),
\end{gather*}
where $M_j(x,\lambda)=C_jT_j(x,\lambda)$ and
\begin{gather*}
C_j=\frac{(b_1/a_j,b_2/a_j,\ldots,b_s/a_j;q)_\infty}{(b_1,b_2,\ldots,b_s;q)_\infty}\prod_{1\le i\le r,\, i\neq j}\frac{(a_i;q)_\infty}{(a_i/a_j;q)_\infty},\\
T_j(x,\lambda)=\frac{\theta_{p}\big(pa_j^kx/\lambda\big)}{\theta_{p}(px/\lambda)}\frac{\theta_q\big((-1)^{1-k}a_j\lambda\big)}{\theta_q\big((-1)^{1-k}\lambda\big)}\frac{\theta_q(-x)}{\theta_q(-a_jx)}\label{T_j}.
\end{gather*}
The coefficients $M_j$ are called $q$-Stokes coefficients. Unlike the case of differential equation, the values of coefficients $M_j$ change continuously dependent on $\lambda\in\mathbb{C}^*\setminus\big[(-1)^k;q\big]$.

Next, we consider to take the limit $q\to1$ in above result. To consider taking the limit of~$q$, we restrict~$q$ to within a real number, and $0<q<1$. When~$q$ is a complex number, the situation is more complicated. For detail, see Sauloy \cite{Sauloy2000}.

\begin{Theorem}\label{Theorem:Main2}We assume that $q\in\mathbb{R}$ and $0<q<1$. Let $\bm{\alpha}=(\alpha_1,\alpha_2,\ldots,\alpha_r)\in\mathbb{C}^r$, $\bm{\beta}=(\beta_1,\beta_2,\ldots,\beta_s)\in\mathbb{C}^s$ and $\alpha_i-\alpha_j\notin\mathbb{Z}$, $i\neq j$. Take a parameter $\lambda\in\mathbb{C}^*\setminus (-1)^k\mathbb{R}_+$. Then for any $x\in\mathbb{C}^*\setminus(-1)^{k+1}\lambda\mathbb{R}_+$, the following equality holds:
\begin{gather}
\lim_{q\to1}{}_rf_s \left(q^{\bm{\alpha}};q^{\bm{\beta}};\lambda;q,\frac{(-1)^kx}{(1-q)^k}\right)\nonumber\\
\qquad {} =C_{\bm{\alpha\beta}}\sum_{j=1}^rC_{\bm{\alpha\beta}}(j)(-x)^{-\alpha_j}{}_{s+1}F_{r-1}\left(\!
\begin{matrix}
\alpha_j,1+\alpha_j-\bm{\beta}\\
1+\alpha_j-\widehat{\bm{\alpha}_j}
\end{matrix}
;\frac{(-1)^k}{x}
\right),\label{Main2}
\end{gather}
where $q^{\bm{\alpha}}=\big(q^{\alpha_1},q^{\alpha_2},\ldots,q^{\alpha_r}\big)$, $q^{\bm{\beta}}=\big(q^{\beta_1},q^{\beta_2},\ldots,q^{\beta_s}\big)$ and $\big|\arg{(-x)}-\arg{(-1)^{k-1}}\lambda\big|<\pi$. $C_{\bm{\alpha\beta}}$~and~$C_{\bm{\alpha\beta}}(j)$ are given by \eqref{Ichinobe_result2}.
\end{Theorem}
We remark that the left hand side of \eqref{Main2} formally converges to ${}_{r}F_{s}(\bm{\alpha};\bm{\beta};x)$ as $q\to1$. Therefore, Theorem~\ref{Theorem:Main2} can be seen as a $q$-analog of Theorem~\ref{Theorem:Ichinobe}.

\section{Proof of Theorem \ref{Theorem:Main1}}\label{Section:Proof_Main1}
We show that $\hat{f}(x)={}_r\varphi_s(\bm{a};\bm{b};q,x)$ is $[\lambda;p]$-summable for any $\lambda\in\mathbb{C}\setminus\big[(-1)^k;q\big]$. Let $g(\xi)$ be the $p$-Borel transform of $\hat{f}(x)$
\begin{gather*}
g(\xi)=(\hat{\mathcal{B}}_{p;1}\hat{f})(\xi)=\sum_{n\ge0}\frac{(a_1,a_2,\ldots,a_r;q)_n}{(b_1,\ldots,b_s;q)_n(q;q)_n}\big\{(-1)^nq^{\frac{n(n-1)}{2}}\big\}^{1+s-r}p^{\frac{n(n-1)}{2}}\xi^n.
\end{gather*}
Then, $g(\xi)$ is convergent in $|\xi|<1$ and again can be expressed by basic hypergeometric series:
\begin{gather*}
g(\xi)=\sum_{n\ge0}\frac{(a_1,a_2,\ldots,a_r;q)_n}{(b_1,\ldots,b_s;q)_n(q;q)_n}\big((-1)^{-k}\xi\big)^n={}_r\varphi_{r-1}\left(
\begin{matrix}
a_1,a_2,\ldots,a_r\\
b_1,b_2,\ldots,b_s,0,\ldots,0
\end{matrix}
;q;(-1)^{-k}\xi
\right).
\end{gather*}
To have the analytic continuation of $g(\xi)$, we use the following proposition.
\begin{Proposition}[Slater \cite{Slater1952} or Gaspar--Rahman {\cite[Section~4.5]{Gasper&Rahman2004}}] The analytic continuation of the basic hypergeometric series ${}_r\varphi_{r-1}(a_1,a_2,\ldots,a_r;b_1,\ldots,b_{r-1};q,\xi)$ is given by
\begin{gather}
{}_r\varphi_{r-1}\left(
\begin{matrix}a_1,\ldots,a_r\\
b_1,\ldots,b_{r-1}
\end{matrix};q,\xi\right)=\frac{(a_2,\ldots,a_r,b_1/a_1,\ldots,b_{r-1}/a_1;q)_\infty}{(b_1,\ldots,b_{r-1},a_2/a_1,\ldots,a_r/a_1;q)_\infty}\frac{\theta_q(-a_1\xi)}{\theta_q(-\xi)}\nonumber\\
\hphantom{{}_r\varphi_{r-1}\left(
\begin{matrix}a_1,\ldots,a_r\\
b_1,\ldots,b_{r-1}
\end{matrix};q,\xi\right)=}{} \times{}_r\varphi_{r-1}\left(
\begin{matrix}a_1,a_1q/b_1,\ldots,a_1q/b_{r-1}\\
a_1q/a_2,\ldots,a_1q/a_r
\end{matrix};q,\frac{qb_1\cdots b_{s}}{a_1\cdots a_r\xi }\right)\nonumber\\
\hphantom{{}_r\varphi_{r-1}\left(
\begin{matrix}a_1,\ldots,a_r\\
b_1,\ldots,b_{r-1}
\end{matrix};q,\xi\right)=}{} +\idem(a_1;a_2,\ldots,a_r).\label{r,r-1AC}
\end{gather}
\end{Proposition}
Letting $b_{s+1},\ldots,b_{r-1}\to 0$ in \eqref{r,r-1AC}, it is seen that $g(\xi)$ is continued analytically to $\xi\in\mathbb{C}^*\setminus\big[(-1)^{k};q\big]$ as follows:
\begin{gather}
g(\xi)=\frac{(a_2,\ldots,a_r,b_1/a_1,\ldots,b_{s}/a_1;q)_\infty}{(b_1\ldots,b_s,a_2/a_1,\ldots a_r/a_1;q)_\infty}\frac{\theta_q((-1)^{k-1}a_1\xi)}{\theta_q((-1)^{k-1}\xi)}\nonumber\\
\hphantom{g(\xi)=}{} \times{}_{s+1}\varphi_{r-1}\left(\tilde{\bm{a}};\tilde{\bm{b}};q,(-1)^k\frac{q^{k+1}b_1\cdots b_s}{a_1^{1-k}a_2\cdots a_r\xi}\right)\nonumber\\
\hphantom{g(\xi)=}{} +\idem(a_1;a_2,\ldots,a_r),\label{ac_borelim}
\end{gather}
where $\tilde{\bm{a}}=(a_1,a_1q/b_1,\ldots,a_1q/b_s)\in\mathbb{C}^{s+1}$ and $\tilde{\bm{b}}=(a_1q/a_2,\ldots,a_1q/a_r)\in\mathbb{C}^{r-1}$. We note that the basic hypergeometric series ${}_{s+1}\varphi_{r-1}$ in~\eqref{ac_borelim} is holomorphic on $\mathbb{C}^*$. In order to obtain \eqref{ac_borelim}, we use
\begin{gather*}
\lim_{b\to0}(aq/b;q)_nb^n=\lim_{b\to0}(b-aq)\big(b-aq^2\big)\cdots\big(b-aq^n\big)=(-a)^nq^{\frac{n(n+1)}{2}}.
\end{gather*}

Now, we shall show that $\big(\hat{\mathcal{B}}_{p;1}\hat{f}\big)(\xi)=g(\xi)\in\mathbb{H}_{p;1}^{[\lambda;p]}$ for any $\lambda\in\mathbb{C}^*\setminus\big[(-1)^k;q\big]$. The following lemma holds.
\begin{Lemma}\label{Lemma:Main1}
The function $g(\xi)$ has a $p$-exponential growth at infinity in
\begin{gather*}
\Omega_\delta:=\mathbb{C}^*\setminus\bigcup_{m\in\mathbb{Z}}\big\{\xi\in\mathbb{C}^* ;\big\vert\xi-(-1)^{k}q^m\vert<\delta\big\vert q^m\vert\big\}
\end{gather*}
for any $\delta>0$.
\end{Lemma}
A proof of this lemma will be given in Section~\ref{Section:Proof_Lemma:Main1}. For any $\lambda\in\mathbb{C}^*\setminus\big[(-1)^k;q\big]$, there exists positive constants $\delta$ and $\varepsilon$ such that
\begin{gather*}
\bigcup_{m\in\mathbb{Z}}\big\{\xi\in\mathbb{C}^* ;\,\vert\xi-\lambda p^m\vert<\varepsilon\vert p^m\lambda\vert\big\}\subset\Omega_\delta.
\end{gather*}
Therefore we obtain $g(\xi)\in\mathbb{H}_{p;1}^{[\lambda;p]}$ for any $\lambda\in\mathbb{C}^*\setminus\big[(-1)^k;q\big]$. In the other words, $\hat{f}(x)={}_r\varphi_s(\bm{a};\bm{b};q,x)$ is $[\lambda;p]$-summable for any $\lambda\in\mathbb{C}^*\setminus\big[(-1)^k;q\big]$.

We consider the $p$-Laplace transform of $g(\xi)$. For any $\lambda\in\mathbb{C}^*\setminus\big[(-1)^k;q\big]$, let ${}_rf_s(\bm{a};\bm{b};\lambda;q,x)$ be the $p$-Laplace transform of $g(\xi)$
\begin{gather*}
{}_rf_s(\bm{a};\bm{b};\lambda;q,x)=\big(\mathcal{L}_{p;1}^{[\lambda;p]}g\big)(x)=\sum_{m\in\mathbb{Z}}\frac{g\big(\lambda p^m\big)}{\theta_p\big(\frac{\lambda p^m}{x}\big)}.
\end{gather*}
Since
\begin{gather*}
\frac{\theta_q\big((-1)^{k-1}a_1\lambda p^m\big)}{\theta_q\big((-1)^{k-1}\lambda p^m\big)}=\frac{\theta_q\big((-1)^{k-1}a_1\lambda q^{km}\big)}{\theta_q\big((-1)^{k-1}\lambda q^{km}\big)}=\frac{1}{a_1^{km}}\frac{\theta_q\big((-1)^{k-1}a_1\lambda\big)}{\theta_q\big((-1)^{k-1}\lambda\big)},
\\
{}_{s+1}\varphi_{r-1}\left(\tilde{\bm{a}};\tilde{\bm{b}};q,(-1)^k\frac{q^{k+1}b_1\cdots b_s}{a_1^{1-k}a_2\cdots a_r\lambda p^m}\right)=\sum_{n\ge0}\frac{(\tilde{\bm{a}};q)_n}{(\tilde{\bm{b}},q;q)_n}p^{\frac{n(n-1)}{2}}\left(\frac{pqb_1\cdots b_s}{a_1^{1-k}a_2\cdots a_r\lambda p^m}\right)^n,
\end{gather*}
and
\begin{gather*}
\theta_p\left(\frac{\lambda p^m}{x}\right)=\left(\frac{\lambda}{x}\right)^{-m}p^{-\frac{m(m-1)}{2}}\theta_p\left(\frac{\lambda}{x}\right),
\end{gather*}
we have
\begin{gather*}
{}_rf_s(\bm{a};\bm{b};\lambda;q,x)= \frac{(a_2,\ldots,a_r,b_1/a_1,\ldots,b_{s}/a_1;q)_\infty}{(b_1\ldots,b_s,a_2/a_1,\ldots, a_r/a_1;q)_\infty}\frac{\theta_q\big((-1)^{k-1}a_1\lambda\big)}{\theta_q\big((-1)^{k-1}\lambda\big)}\frac{1}{\theta_p(\lambda/x)}\\
\hphantom{{}_rf_s(\bm{a};\bm{b};\lambda;q,x)=}{} \times\sum_{n\ge0}\frac{(\tilde{\bm{a}};q)_n}{(\tilde{\bm{b}},q;q)_n}\left(\frac{qb_1\cdots b_s}{a_1^{1-k}a_2\cdots a_r\lambda}\right)^n \sum_{m\in\mathbb{Z}}\left(\frac{\lambda}{a_1^kx}\right)^mp^{\frac{n(n-1)}{2}+n-mn+\frac{m(m-1)}{2}}\\
\hphantom{{}_rf_s(\bm{a};\bm{b};\lambda;q,x)=}{} +\idem(a_1;a_2,\ldots,a_r).
\end{gather*}
Here we remark that
\begin{gather*}
\sum_{m\in\mathbb{Z}}\left(\frac{\lambda}{a_1^kx}\right)^mp^{\frac{n(n-1)}{2}+n-mn+\frac{m(m-1)}{2}} =\left(\frac{\lambda}{a_1^kx}\right)^n\sum_{m\in\mathbb{Z}}p^{\frac{(m-n)(m-n-1)}{2}}\left(\frac{\lambda}{a_1^kx}\right)^{m-n}\\
\hphantom{\sum_{m\in\mathbb{Z}}\left(\frac{\lambda}{a_1^kx}\right)^mp^{\frac{n(n-1)}{2}+n-mn+\frac{m(m-1)}{2}} }{}
=\left(\frac{\lambda}{a_1^kx}\right)^n\theta_p\left(\frac{\lambda}{a_1^kx}\right).
\end{gather*}
Therefore we obtain
\begin{gather*}
{}_rf_s(\bm{a};\bm{b};\lambda;q,x)= \frac{(a_2,\ldots,a_r,b_1/a_1,\ldots,b_{s}/a_1;q)_\infty}{(b_1\ldots,b_s,a_2/a_1,\ldots a_r/a_1;q)_\infty}\frac{\theta_q\big((-1)^{k-1}a_1\lambda\big)}{\theta_q\big((-1)^{k-1}\lambda\big)}\frac{\theta_p\big(\lambda/a_1^kx\big)}{\theta_p(\lambda/x)}\\
\hphantom{{}_rf_s(\bm{a};\bm{b};\lambda;q,x)=}{} \times\sum_{n\ge0}\frac{(\tilde{\bm{a}};q)_n}{(\tilde{\bm{b}},q;q)_n}\left(\frac{qb_1\cdots b_s}{a_1a_2\cdots a_rx}\right)^n\\
\hphantom{{}_rf_s(\bm{a};\bm{b};\lambda;q,x)=}{} +\idem(a_1;a_2,\ldots,a_r).
\end{gather*}
From
\begin{gather*}
\frac{\theta_p\big(\lambda/a_1^kx\big)}{\theta_p(\lambda/x)}=\frac{\theta_p\big(pa_1^kx/\lambda\big)}{\theta_p(px/\lambda)}
\end{gather*}
(is seen from equality \eqref{theta_inverse}) and
\begin{gather*}
\sum_{n\ge0}\frac{(\tilde{\bm{a}};q)_n}{(\tilde{\bm{b}},q;q)_n}\left(\frac{qb_1\cdots b_s}{a_1a_2\cdots a_rx}\right)^n={}_r\varphi_{r-1}\left(
\begin{matrix}\tilde{\bm{a}},\bm{0}_k\\
\tilde{\bm{b}}
\end{matrix};q,\frac{qb_1\cdots b_s}{a_1a_2\cdots a_rx}\right),
\end{gather*}
we finish the proof.
\section{Proof of Lemma \ref{Lemma:Main1}}\label{Section:Proof_Lemma:Main1}
In this section, we give a proof of Lemma \ref{Lemma:Main1}, that is, we show
\begin{gather*}
\vert g(\xi)\vert\le L \theta_{\vert p\vert}(M\vert \xi\vert),\qquad \xi\in\Omega_\delta,
\end{gather*}
where $L$ and $M$ are positive constants.

As easily seen, the following inequality holds.
\begin{gather*}
\theta_{|p|}(A|\xi|)+\theta_{|p|}(B|\xi|)\le \theta_{|p|}(C|\xi|),
\end{gather*}
where $A, B>0$ and $C=A+B$. Hence, showing that the first term of \eqref{ac_borelim} has $p$-exponential growth at infinity in $\Omega_\delta$ is sufficient to prove Lemma \ref{Lemma:Main1}. We write the first term of \eqref{ac_borelim} as~$g_1(\xi)$, i.e.,
\begin{gather*}
g_1(\xi):=\frac{(a_2,\ldots,a_r,b_1/a_1,\ldots,b_{s}/a_1;q)_\infty}{(b_1\ldots,b_s,a_2/a_1,\ldots a_r/a_1;q)_\infty}\frac{\theta_q\big((-1)^{k-1}a_1\xi\big)}{\theta_q\big((-1)^{k-1}\xi\big)}\\
\hphantom{g_1(\xi):=}{}\times {}_{s+1}\varphi_{r-1}\left(\tilde{\bm{a}};\tilde{\bm{b}};q;(-1)^k\frac{q^{k+1}b_1\cdots b_s}{a_1^{1-k}a_2\cdots a_r\xi}\right),
\end{gather*}
where $\tilde{\bm{a}}=(a_1,a_1q/b_1,\ldots,a_1q/b_s)\in\mathbb{C}^{s+1}$ and $\tilde{\bm{b}}=(a_1q/a_2,\ldots,a_1q/a_r)\in\mathbb{C}^{r-1}$.

We show that basic hypergeometric series ${}_{s+1}\varphi_{r-1}$ has $p$-exponential growth. It holds that
\begin{gather*}
|(\alpha;q)_n| =|1-\alpha||1-\alpha q|\cdots\big|1-\alpha q^{n-1}\big|\le(1+|\alpha|)^n,\\
|(a_iq/a_j;q)_n| =|1-a_iq/a_j|\big|1-a_iq^2/a_j\big|\cdots\big|1-a_iq^n/a_j\big|\ge \left(\inf_{k\in\mathbb{N}^*}\big\{\big|1-a_iq^k/a_j\big|\big\}\right)^n>0
\end{gather*}
for any $\alpha\in\mathbb{C}$ and $1\le i,j\le r$. Therefore we obtain
\begin{align*}
\left|{}_{s+1}\varphi_{r-1}\left(\tilde{\bm{a}};\tilde{\bm{b}};q;(-1)^k\frac{q^{k+1}b_1\cdots b_s}{a_1^{1-k}a_2\cdots a_r\xi}\right)\right|&\le\sum_{n\ge0}\frac{|(\tilde{\bm{a}};q)_n|}{|(\tilde{\bm{b}},q;q)_n|}|p|^{\frac{n(n+1)}{2}}\left|\frac{qb_1\cdots b_s}{a_1^{1-k}a_2\cdots a_r\xi}\right|^n\\
&\le\sum_{n\ge0}|p|^{\frac{n(n+1)}{2}}\left(\frac{1}{M|\xi|}\right)^n\\
&\le\sum_{n\in\mathbb{Z}}|p|^{\frac{n(n+1)}{2}}\left(\frac{1}{M|\xi|}\right)^n=\theta_{|p|}(M|\xi|),
\end{align*}
where $M$ is a positive constant. Hence we have
\begin{gather}\label{g1_estimate_1}
|g_1(\xi)|\le C\left|\frac{\theta_q\big((-1)^{k-1}a_1\xi\big)}{\theta_q\big((-1)^{k-1}\xi\big)}\right|\theta_{|p|}(M|\xi|).
\end{gather}
Here $C$ is a positive constant.

Now, for any $\xi \in\Omega_\delta$, there exist an integer $n$ and
\begin{gather*}
\xi_0\in\Omega_\delta \cap\{\xi\in\mathbb{C}^* ;\,|p|\le|\xi|\le1\}
\end{gather*}
such that $\xi=p^n\xi_0$.
Substituting $\xi=p^n\xi_0$ into \eqref{g1_estimate_1}, we have
\begin{gather*}
|g_1(\xi)|\le C\left|\frac{\theta_q((-1)^{k-1}a_1\xi_0)}{\theta_q((-1)^{k-1}\xi_0)}\right|\big(M|a_1|^k\big)^{-n}|p|^{-\frac{n(n-1)}{2}}\theta_{|p|}(M|\xi_0|).
\end{gather*}
Since
\begin{gather*}
\big(M|a_1|^k\big)^{-n}|p|^{-\frac{n(n-1)}{2}}=\frac{\theta_{|p|}\big(M|a_1|^k|p|^n|\xi_0|\big)}{\theta_{|p|}\big(M|a_1|^k|\xi_0|\big)}=\frac{\theta_{|p|}\big(M|a_1|^k|\xi|\big)}{\theta_{|p|}\big(M|a_1|^k|\xi_0|\big)}
\end{gather*}
holds, we obtain
\begin{gather*}
|g_1(\xi)|\le C\left|\frac{\theta_q\big((-1)^{k-1}a_1\xi_0\big)}{\theta_q\big((-1)^{k-1}\xi_0\big)}\right|\frac{\theta_{|p|}(M|\xi_0|)}{\theta_{|p|}\big(M|a_1|^k|\xi_0|\big)}\theta_{|p|}\big(M|a_1|^k|\xi|\big).
\end{gather*}
Now, the function of $\xi_0$
\begin{gather*}
\left|\frac{\theta_q\big((-1)^{k-1}a_1\xi_0\big)}{\theta_q\big((-1)^{k-1}\xi_0\big)}\right|\frac{\theta_{|p|}(M|\xi_0|)}{\theta_{|p|}\big(M|a_1|^k|\xi_0|\big)}
\end{gather*}
is holomorphic on the compact set $\Omega_\delta \cap\{\xi\in\mathbb{C}^* ;\,|p|\le|\xi|\le1\}$. Hence there exists a positive constant $L$ such that
\begin{gather*}
\max_{\Omega_\delta \cap\{\xi\in\mathbb{C}^* ;\,|p|\le|\xi|\le1\}}\left|\frac{\theta_q\big((-1)^{k-1}a_1\xi_0\big)}{\theta_q\big((-1)^{k-1}\xi_0\big)}\right|\frac{\theta_{|p|}(M|\xi_0|)}{\theta_{|p|}\big(M|a_1|^k|\xi_0|\big)}\le L.
\end{gather*}

\section{Proof of Theorem \ref{Theorem:Main2}}\label{Section:Proof_Main2}
A proof of Theorem \ref{Theorem:Main2} is obtained from the following proposition.
\begin{Proposition}[Askey {\cite[Section~5]{Askey1978}}]
We assume that $q\in\mathbb{R}$ and $0<q<1$. For any $x\in\mathbb{C}^*\setminus\mathbb{R}_-$, we have
\begin{gather}\label{theta_reduction1}
\lim_{q\to1}\frac{\theta_q\big(q^\beta x\big)}{\theta_q(q^\alpha x)}=x^{\alpha-\beta}
\end{gather}
and
\begin{gather}\label{theta_reduction2}
\lim_{q\to1}\frac{\theta_p\big(\frac{p^\alpha x}{(1-q)^k}\big)}{\theta_p\big(\frac{p^\beta x}{(1-q)^k}\big)}(1-q)^{k(\beta-\alpha)}=x^{\beta-\alpha}.
\end{gather}
Here we assume $|\arg{x}|<\pi$.
\end{Proposition}
Let us give a proof of Theorem~\ref{Theorem:Main2}. Substituting $\bm{a}=q^{\bm{\alpha}}$, $\bm{b}=q^{\bm{\beta}}$ and $x=(-1)^kx/(1-q)^k$ into~\eqref{Main1}, we have
\begin{gather*}
{}_rf_s\left(q^{\bm{\alpha}};q^{\bm{\beta}};\lambda;q,\frac{(-1)^kx}{(1-q)^k}\right) =\frac{\Gamma_q(\beta_1)\cdots\Gamma_q(\beta_s)\Gamma_q(\alpha_2-\alpha_1)\cdots\Gamma_q(\alpha_r-\alpha_1)}{\Gamma_q(\alpha_2)\cdots\Gamma_q(\alpha_r)\Gamma_q(\beta_1-\alpha_1)\cdots\Gamma_q(\beta_s-\alpha_1)}(1-q)^{-k\alpha_1}\\
\hphantom{{}_rf_s\left(q^{\bm{\alpha}};q^{\bm{\beta}};\lambda;q,\frac{(-1)^kx}{(1-q)^k}\right) =}{} \times\frac{\theta_{p}\left(\frac{p^{\alpha_1+1}(-1)^kx}{\lambda(1-q)^k}\right)}{\theta_{p}\left(\frac{p(-1)^kx}{\lambda(1-q)^k}\right)}\frac{\theta_q((-1)^{1-k}q^{\alpha_1}\lambda)}{\theta_q((-1)^{1-k}\lambda)}\\
\hphantom{{}_rf_s\left(q^{\bm{\alpha}};q^{\bm{\beta}};\lambda;q,\frac{(-1)^kx}{(1-q)^k}\right) =}{} \times{}_r\varphi_{r-1}\left(\begin{matrix}q^{\alpha_1},q^{1+\alpha_1-\beta_1},\ldots,q^{1+\alpha_1-\beta_s},\bm{0}_k\\q^{1+\alpha_1-\alpha_2},\ldots,q^{1+\alpha_1-\alpha_r}
\end{matrix} ;q;q^{\gamma}\frac{(1-q)^k}{(-1)^kx}\right)\\
\hphantom{{}_rf_s\left(q^{\bm{\alpha}};q^{\bm{\beta}};\lambda;q,\frac{(-1)^kx}{(1-q)^k}\right) =}{} +\idem\big(q^{\alpha_1};q^{\alpha_2},\ldots,q^{\alpha_r}\big),
\end{gather*}
where $\gamma=1+\beta_1+\cdots+\beta_s-(\alpha_1+\cdots+\alpha_r)$. Since we have
\begin{gather*}
\lim_{q\to1}\frac{\Gamma_q(\beta_1)\cdots\Gamma_q(\beta_s)\Gamma_q(\alpha_2-\alpha_1)\cdots\Gamma_q(\alpha_r-\alpha_1)}{\Gamma_q(\alpha_2)\cdots\Gamma_q(\alpha_r)\Gamma_q(\beta_1-\alpha_1)\cdots\Gamma_q(\beta_s-\alpha_1)}\\
\qquad{} =\frac{\Gamma(\beta_1)\cdots\Gamma(\beta_s)\Gamma(\alpha_2-\alpha_1)\cdots\Gamma(\alpha_r-\alpha_1)}{\Gamma(\alpha_2)\cdots\Gamma(\alpha_r)\Gamma(\beta_1-\alpha_1)\cdots\Gamma(\beta_s-\alpha_1)},\\
\lim_{q\to1}\frac{\theta_q\big((-1)^{1-k}q^{\alpha_1}\lambda\big)}{\theta_q\big((-1)^{1-k}\lambda\big)}=\big\{(-1)^{1-k}\lambda\big\}^{-\alpha_1},\\
\lim_{q\to1}\frac{\theta_{p}\big(\frac{p^{\alpha_1+1}(-1)^kx}{\lambda(1-q)^k}\big)}{\theta_{p}\big(\frac{p(-1)^kx}{\lambda(1-q)^k}\big)}(1-q)^{-k\alpha_1}=\left\{(-1)^k\frac{x}{\lambda}\right\}^{-\alpha_1},\\
\lim_{q\to1}{}_r\varphi_{r-1}\left(\begin{matrix}q^{\alpha_1},q^{1+\alpha_1-\beta_1},\ldots,q^{1+\alpha_1-\beta_s},\bm{0}_k\\q^{1+\alpha_1-\alpha_2},\ldots,q^{1+\alpha_1-\alpha_r}
\end{matrix} ;q;q^{\gamma}\frac{(1-q)^k}{(-1)^kx}\right)\\
\qquad{} ={}_{s+1}F_{r-1}\left(\begin{matrix}\alpha_1,1+\alpha_1-\bm{\beta}\\1+\alpha_1-\widehat{\bm{\alpha}_1}
\end{matrix} ;\frac{(-1)^k}{x}\right)
\end{gather*}
from \eqref{qgamma_reduction}, \eqref{theta_reduction1}, \eqref{theta_reduction2} and
\begin{gather*}
\lim_{q\to1}\frac{(q^\alpha;q)_n}{(1-q)^n}=(\alpha)_n
\end{gather*}
respectively, we obtain \eqref{Main2}.

\subsection*{Acknowledgements}

The author would like to thank the referees for their helpful suggestions and variable comments. Additionally, the author is grateful for that one of the referee let him know the existence of the papers \cite{Dreyfus2015_2} and \cite{Sauloy2000}.

\pdfbookmark[1]{References}{ref}
\LastPageEnding

\end{document}